\documentclass[a4paper]{amsart}

\usepackage{amssymb}
\usepackage[colorlinks=true, linkcolor=blue]{hyperref}
\usepackage[all]{xypic}

\usepackage{amsmath}
\usepackage{amsfonts}
\usepackage{amssymb}
\usepackage[english]{babel}
\usepackage{amsthm}
\usepackage{extarrows}
\usepackage{amscd}
\usepackage{tikz}
\usepackage{mathdots}
\usetikzlibrary{arrows,shapes}

\newif\ifFIRST\newdimen\MAXright\MAXright0pt
\def\dynkin{\mbox\bgroup\scriptsize
\let\SMALL\relax
\FIRSTtrue\hskip.5em%
\setbox1\hbox{$\diagup$}%
\setbox2\hbox{$\diagdown$}%
\setbox0\hbox to2\wd1{\hrulefill}%
\setbox8\hbox to.5\wd1{\hrulefill}%
\setbox3\hbox{$\bullet$}%
\setbox4\hbox{$\times$}%
\def\root##1{\ifFIRST\setbox5\hbox{$##1$}\ifdim\wd5>1.3em%
\hskip-.5em\hskip.5\wd5\fi\fi\FIRSTfalse%
\hskip-.25em\raise1.5\wd3\hbox to0pt{\hss\hskip.45em$%
##1$\hss}\copy3\hskip-.25em\setbox6\hbox{$##1$}%
\MAXright\wd6}%
\def\droot##1{\ifFIRST\setbox5\hbox{$##1$}\ifdim\wd5>1.3em%
\hskip-.5em\hskip.5\wd5\fi\fi\FIRSTfalse%
\hskip-.25em\lower1.8\wd3\hbox to0pt{\hss\hskip.45em$%
##1$\hss}\copy3\hskip-.25em\setbox6\hbox{$##1$}%
\MAXright\wd6}%
\def\rroot##1{\hskip-.25em\copy3\hbox to0pt{\hskip.3em$##1$\hss}%
\hskip-.25em\setbox6\hbox{\hskip.6em$##1##1$}%
\MAXright\wd6}%
\def\norroot##1{\hskip-.36em\copy4\hbox to0pt{\hskip.3em$##1$\hss}%
\hskip-.48em\setbox6\hbox{\hskip.6em$##1##1$}%
\MAXright\wd6}%
\def\noroot##1{\ifFIRST\setbox5\hbox{$##1$}\ifdim\wd5>1.3em%
\hskip-.5em\hskip.5\wd5\fi\fi\FIRSTfalse%
\hskip-.36em\raise1.5\wd3\hbox to0pt{\hss\hskip.6em$%
##1$\hss}\copy4\hskip-.38em\setbox6\hbox{$##1$}%
\MAXright\wd6}%
\def\nodroot##1{\ifFIRST\setbox5\hbox{$##1$}\ifdim\wd5>1.3em%
\hskip-.5em\hskip.5\wd5\fi\fi\FIRSTfalse%
\hskip-.36em\lower1.8\wd3\hbox to0pt{\hss\hskip.6em$%
##1$\hss}\copy4\hskip-.38em\setbox6\hbox{$##1$}%
\MAXright\wd6}%
\def\link{\raise.25em\copy0}%
\def\minilink{\raise.25em\copy8}%
\def\llink##1{\raise.35em\copy0\hskip-\wd0%
\raise.12em\copy0\hskip-.5\wd0\hbox to0pt{\hss$##1$\hss}\hskip.5\wd0}%
\def\lllink##1{\raise.40em\copy0\hskip-\wd0\raise.25em\copy0\hskip-\wd0%
\raise.10em\copy0\hskip-.5\wd0\hbox to0pt{\hss$##1$\hss}\hskip.5\wd0}%
\def\rootupright##1{\hbox to0pt{\raise.45em\copy1\hskip-.25em\raise1.3\ht1%
\hbox{\copy3\hskip.3em$##1$}\hss}%
\setbox6\hbox{\hskip.6em\copy1\copy1$##1##1$}%
\ifdim\MAXright<\wd6\MAXright\wd6\fi}%
\def\norootupright##1{\hbox to0pt{\raise.45em\copy1\hskip-.25em\raise1.3\ht1%
\hbox{\kern-.1em\copy4\hskip.3em$##1$}\hss}%
\setbox6\hbox{\hskip.6em\copy1\copy1$##1##1$}%
\ifdim\MAXright<\wd6\MAXright\wd6\fi}%
\def\rootdownright##1{\hbox to0pt{\raise-.5em\copy2\hskip-.25em\raise-1.35\ht1%
\hbox{\copy3\hskip.3em$##1$}\hss}\setbox6%
\hbox{\hskip.6em\copy2\copy2$##1##1$}%
\ifdim\MAXright<\wd6\MAXright\wd6\fi}%
\def\norootdownright##1{\hbox to0pt{\raise-.5em\copy2\hskip-.25em\raise-1.35\ht1%
\hbox{\kern-.1em\copy4\hskip.3em$##1$}\hss}\setbox6%
\hbox{\hskip.6em\copy2\copy2$##1##1$}%
\ifdim\MAXright<\wd6\MAXright\wd6\fi}%
\def\rootdown##1{\hbox to0pt{\hskip-.05em\vrule height.25em depth.65em%
\hskip-.25em\raise-.95em\hbox{\copy3\hskip.3em$##1$}\hss}%
\setbox6\hbox{$##1$}%
\ifdim\MAXright<\wd6\MAXright\wd6\fi}%
\def\norootdown##1{\hbox to0pt{\hskip-.05em\vrule height.25em depth.65em%
\hskip-.25em\raise-.95em\hbox{\copy4\hskip.3em$##1$}\hss}%
\setbox6\hbox{$##1$}%
\ifdim\MAXright<\wd6\MAXright\wd6\fi}%
\def\dots{\raise.25em\copy8\hskip.25em\raisebox{-0.025em}{$\cdots$}\hskip.10em%
\raise.25em\copy8}}%
\def\enddynkin{\ifdim\MAXright>1em\hskip.5\MAXright\else\hskip.5em\fi\egroup}%
 
This macro is used for creation of Dynkin diagrams in AMSTeX.
The usage:
In math-mode, between \dynkin and \enddynkin,
\root#1 makes a node with label #1 over this node,
\rroot#1 and \droot do the same with labels on right and down,
\noroot and \norroot do the same with cross instead of bullet,
\link links two adjacent \root's, \llink#1 links two adjacent nodes
with double-line with #1 (probably '>' or '<') in the middle, \lllink is does
same with three lines, \dots does \cdots between the two adjacent nodes,
\rootdown#1 makes a linked node down labeled by #1, \rootupright#1 and
\rootdownright#1 are labeled linked nodes in the indicated directions 
which do not count for the horizontal dimensions (so the user has to adjust
ad hoc the break after the diagram).

Examples:
$\dynkin \root{a_1}\link\root{a_2}\dots\root{a_{n-1}}\link\root{a_n}
\enddynkin$
or
$\dynkin \root{}\lllink>\root{}\enddynkin$
or
$\dynkin \root{a}\link\root{b}\rootupright{c}\rootdownright{d}\enddynkin$
or
$\dynkin \root{}\link\root{}\rootupright{}\link\root{}\enddynkin$.

\newcommand{\dyn}[5]{\dynkin
\root{#1}\link\root{#2}\link\noroot{#3}\link\root{#4}\link\root{#5}
\enddynkin}

\newcommand{\wgt}[2]{{(#1|#2)}}

\newcommand{\bR}{\mathbb{R}}

\newcommand{\bN}{\mathbb{N}}
\newcommand{\bZ}{\mathbb{Z}}
\newcommand{\bE}{\mathbb{E}}

\newcommand{\bV}{\mathbb{V}}

\newcommand{\bW}{\mathbb{W}}
\newcommand{\bF}{\mathbb{F}}

\newcommand{\al}{\alpha}
\newcommand{\be}{\beta}

\newcommand{\mfg}{\mathfrak{g}}
\newcommand{\mfp}{\mathfrak{p}}

\newcommand{\la}{\lambda}

\newcommand{\pa}{\partial}

\newtheorem{theorem}{Theorem}
\newtheorem{lemma}{Lemma}
\newtheorem{proposition}{Proposition}

\begin{document}
	
\title[Strongly invariant differential operators on $Gr(3,3)$]{Strongly 
invariant differential operators on parabolic geometries modeled on $Gr(3,3)$}

\author{Jan Slov\'ak and Vladim\'{\i}r Sou\v cek}
\address{J.S.: Department of Mathematics and Statistics, Faculty of Science,
Masaryk University, Kotl\'a\v rsk\'a 2a, 611~37~Brno, Czech Republic
\newline\indent V.S.: Mathematical Institute of Charles University,
Sokolovská 83, 186 75 Praha, Czech Republic}

\email{slovak@math.muni.cz, soucek@karlin.mff.cuni.cz}

%

\thanks{Both authors
acknowledge the support by the grant GX19-28628X of the Czech Grant Agency.
The first author was also supported by the projects CaLIGOLA and CaLiForNIA, Horizon Europe, 
project id 101086123 and 101119552. The second author acknowledges the support by the grant 26-20231L financed by the Czech Grant Agency. The first author would also like to thank the Isaac Newton
Institute for Mathematical Sciences for the support and hospitality during
the program Twistor Theory.}


	
\begin{abstract}
We consider the curved geometries modeled on the homogeneous space
$G/P,$ where $G=SL(6,\bR)$ acts transitively on the Grassmannian $Gr(3,3)$
of three-dimensional subspaces in $\bR^6$, and $P$ is the corresponding
isotropic subgroup. We classify the strongly
invariant operators between sections of vector bundles induced on such
geometries by
irreducible $P$-modules, i.e., those obtained via 
homomorphisms of semi-holonomic Verma modules.		
\end{abstract}

\maketitle	


\section{Introduction} The 4-dimensional conformal Riemannian geometry
twistor calculus, going back to Roger Penrose and his group, has got essential in
modern theories of classical gravity. Over decades, it was generalized
in two main directions -- the conformal geometry and
the (almost) quaternionic geometry, in general dimensions. 
The quaternionic geometry is closely related to
the $(2,n)$-Grassmanians, overlapping with the original $(2,2)$-Grassmannian
story in the most classical spinor and twistor calculus.

We suggest to explore the third obvious possibility here --- the
$(n,n)$-Grassmannian geometries.  Compared to the 4-dimensional case, this
choice has got several nice similarities.  First, the tangent space is a
tensor product of an auxiliary $n$-dimensional bundle $\mathcal S$ and its
dual, whose top forms are identified, i.e., $\Lambda^n \mathcal S \simeq
\Lambda^n\mathcal S^*$ as line bundles.  Second, there is the canonical
Cartan connection coming from the standard normalization of the parabolic
geometries, with all the relevant calculus at hand.  Finally, there are also
two essential components of the Cartan curvature here, which are both
of the torsion type for all $n\ge3$.

All three series of geometries are belonging to the category of
($|1|$-graded) parabolic geometries, see \cite{CS} for careful exposition. 
Thus, the general calculus applies and it is well known that the main tools
of these calculi are of algebraic and (co)homological character. In
particular, there is the complete understanding of patterns of all 
invariant linear differential operators between the
homogeneous vector bundles for the flat Klein geometry models. They correspond
to the well known patterns of (generalized) Verma module homomorphisms, a
story going back to Gelfand, Kostant and others, many decades back. 

There are several approaches how to extend and understand these patterns for 
the curved Cartan geometries, cf. \cite{CSS, CD, CS07, CS12}. 
The work \cite{ES} provides the purely algebraic picture by introducing the
semi-holonomic induced modules and extending the celebrated
Jantzen-Zuckermann translation principle. 

We aim at classification of strongly invariant operators between sections of
vector bundles induced by irreducible $P$-modules.  In general, the
concept of \emph{strongly} invariant operators is related to the tractor
calculi.  It was introduced in the realm of conformal geometry and invariant
operators on scalar functions which can be directly viewed as operators on
any tractors, see the survey in \cite{E96}.  This always happens if the
operators on the flat model are given by $P$-module homomorphisms on the
semi-holonomic jets covering the relevant $P$-module homomorphisms on jets. 
In this work, we understand the \emph{strongly invariant operators} in this
sense only, and these operators also automatically extend to all curved
geometries.  Clearly, this corresponds to classification of homomorphisms of
semi-holonomic Verma modules, as introduced in \cite{ES}.
  
All these concepts and results are summarized by the authors in the recent survey article
\cite{SlS}.  In this paper, we shall rely on all the facts
explained there. Facts from representation theory could be also collected
from \cite{vogan}.

\section{Basic concepts and notation}

Grassmannian manifolds are homogeneous spaces $Gr(m,n)=G/P,$ where
$G=SL(m+n,\bR)$ and $P$ is its parabolic subgroup consisting of the $m+n$ block-wise upper triangular matrices. The Lie algebra $\mfp$
of $P$ has got $\mathfrak{sl}(m,\bR)\times \mathfrak{sl}(n, \bR)$ as the simple
part of the Levi subalgebra of $\mfp$. 

We shall deal with normal parabolic geometries $M$ modelled on the homogeneous space
$G/P$, where $G=SL(6,\bR)$ acts transitively on the Grassmannian $Gr(3,3)$
of three-dimensional subspaces in $\bR^6$ and $P$ is the corresponding
isotropic subgroup. 
 

Let $\mfg$ be the Lie algebra of $G$. Irreducible representations of $\mfg$
are labeled by highest weights $\la=(\la_1\la_2\la_3\la_4\la_5\la_6)$
with $\la_i\in\bZ$; $i=1,\ldots,6$; $\la_{i}\geq\la_{i+1}$; $i=1,\ldots,5$, modulo the
equivalence relation: $\la\simeq\la'$ if and only if 
$\la'_i=\la_i+c$; $i=1,\ldots,6$, for any fixed integer $c\in \mathbb Z$.

Often we shall fix a representative of the equivalence relation by the
additional condition $\la_6=0$.

The differences $\la_i-\la_{i+1}$ correspond to the coefficients in the
alternative description of the highest weights as 
linear combinations of the fundamental weights. Thus, $\rho=(543210)$ is 
the lowest weight in the interior of the dominant Weyl chamber of
$\mathfrak(g)$. 
The irreducible
$\mfg$-module with the highest weight $\la$ will be denoted by $\bW_\al$, 
where $\al=\la+\rho$.
The components of $\al $ satisfy
$\al_i\in\bZ$, $i=1,\ldots,6$; $\al_{i+1}<\al_i$; $i=1,\ldots,5$.

The Lie algebra $\mfg $ is $|1|$-graded, it decomposes under the adjoint
action of $\mfg_0$ into irreducible components as
$\mfg=\mfg_{-1}\oplus\mfg_0\oplus\mfg_1$. 

Irreducible representations of the parabolic subalgebra $\mfp$ are given by
irreducible modules for the Levi factor $\mfg_0$, 
with the trivial action of the nilpotent radical $\mathfrak g_1$.
They are classified by their highest weights, and we can consider them
as the so-called $\mathfrak p$-dominant weights of $\mathfrak g$. 
We shall again denote such a module with highest $\la$
as $\bV_\al,\al=\la+\rho$, and write them with a vertical bar in the middle. 
Components of such $\al$ satisfy the inequalities:
$$
(\al_1\al_2\al_3|\al_4\al_5\al_6);\;\al_1>\al_2>\al_3;\;\al_4>\al_5>\al_6;\;
\al_i\in\bZ,\;i=1,\ldots,6.
$$
with the same equivalence relation as above. 

Notice that the difference
$\al_3-\al_4$ is responsible for the resulting action of the center of
$\mathfrak g_0$, as follows. 
Elements of the center of $\mfg_0$
act by scalar multiplication on irreducible $\mfp$-modules. There is the unique
element in the center of $\mfg_0$ whose adjoint action $\varphi$ on 
$\mathfrak g$ is 
$\varphi(\mfg_j)=j,j=-1,0,1$.  For irreducible $\mfp$-modules with highest
weight $\la$, the value of the
corresponding scalar is given by
$$
\varphi(\al)=\varphi(\mathbb V_\al)=\frac{1}{2}(\la_1+\la_2+\la_3-\la_4-\la_5-\la_6),
$$
and we have to subtract $\varphi(\rho)=\frac12 9$, when using the 
$\al = \la+\rho$ instead. 

A general finite dimensional $P$-module $\bE$ has a composition series
$\bV=\bV_{\al_0}+\ldots + \bV_{\al_k}$, 
where $\varphi(\bV_{\al_j})=\varphi(\al_0)+j$, $j=1,\ldots, k$.

Given an irreducible $P$-module $\mathbb V_\al$, the (generalized) 
Verma module is defined by 
$V(\mathbb V_\al)=\mathfrak{U}(\mfg)\otimes_{\mathfrak{U}(\mfp)}\mathbb
(\mathbb V_\al)^*$. Shortly, $V(\mathbb{V}_\alpha)=V(\alpha)$.

As well known, the generalized Verma modules $V(\al)$ are the topological
duals to the infinite jets of section of the homogeneous bundles defined by
the $P$-modules $\mathbb V_\al$.  Thus the morphisms between them correspond
to the invariant linear differential operators, see \cite{SlS}.

The \emph{order} of a homomorphism $\Phi: V(\al)\to V(\be)$ is defined as
$\varphi(\be)-\varphi(\al)$ (and this clearly corresponds to the order of
the corresponding differential operators).

The Verma modules are $\mathfrak U(\mathfrak g)$-modules. Clearly,
nontrivial homomorphisms between two Verma modules may exist only if the
actions of the center $\mathfrak Z(\mathfrak g)$ of $\mathfrak U(\mathfrak
g)$ on them coincide.  The action of $\mathfrak Z(\mathfrak g)$ is called
the \emph{infinitesimal character} and it is always shared just by Verma
modules for the weights $\al$ on one orbit of the affine Weyl group action
(and this is the reason why the notation $\al=\la+\rho$ is so useful). The
action of the Weyl group on our weights $\al$ is given by permutations
of its components (another reason for choosing this notation). Thus, the 
infinitesimal characters of $V(\al)$ and $V(\al')$ coincide if and
only if
there is a permutation $\pi$ of the set $\{1,2,3,4,5,6\}$ such that
$$ \al'_i=\al_{\pi(i)},\; i=1,\ldots,6.$$

The infinitesimal character of $V(\al)$ is called \emph{regular} if the affine 
orbit contains only one $\mathfrak g$-dominant weight inside the dominant Weyl
chamber, and this happens if and
only if all components $\al_i$, $i=1,\ldots,6$ are mutually different. 
We call
the infinitesimal character of $V(\al)$ $|1|$-singular (resp. 
$|2|$-singular, resp.  $|3|$-singular) if and only if one pair (resp.  two
pairs, resp.  three pairs) of components coincide (which means that the only
$\mathfrak g$-dominant weight $\al$ in the affine orbit 
is on intersections of 2, 3, or 4 walls of the
dominant Weyl chamber).

For our purposes, we do not need to
consider Verma modules with $|3|$-singular infinitesimal characters, because
there are no nontrivial homomorphisms between such Verma modules.

Another important observation is that for irreducible $\mathfrak p$-modules
$\mathbb F^*$ and $\mathbb E^*$, 
the Verma modules are again highest weight modules.
Moreover, as $\mathfrak
U(\mathfrak p)$-modules they come with a
filtration 
$$
\mathbb R+V_1(\be)+V_2(\be) + \dots + V_k(\be) + \dots 
,$$
where the quotients are equal to $V_k(\be)/V_{k-1}(\be)=S^k\mathfrak
g_{-1}\otimes \mathbb E^*$, as $\mathfrak g_0$-modules. Thus, every
homomorphism $\psi:V(\al)\to V(\be)$ must be completely determined by the
embedding $\mathbb F^* \to V(\be)$ and this must sit in a $V_k(\be)$ with
the smallest possible $k$ for our $|1|$-graded algebra. Of course, $k$
equals the order of the homomorphism.

The filtrations define the short exact sequences 
$$
0\to V_{k-1}(\be)\to V_k(\be)\to S^k\mathfrak g_{-1}\otimes \mathbb
E^* \to 0 
$$ 
and the composition of the embedding of $\mathbb F^*$ with the projection 
defines the \emph{symbol} of the homomorphism $\psi$, 
$\sigma(\psi):\mathbb F^*\to S^k\mathfrak
g_{-1}\otimes \mathbb E^*$. 
Notice that for $|1|$-graded geometries, the homomorphisms are completely
determined by their symbols.
All these facts are explained in more detail in \cite{SlS}.

\section{Classification of homomorphism between Verma modules}
On the way to our aim, the discussion of the homomorphisms of the
semi-holonomic Verma modules, we need the well known 
complete description of all
homomorphisms between the classical Verma modules, i.e., the complete
classification of all linear invariant operators between the homogenous
bundles on our Grassmannian Klein model $Gr(3,3)$.

Let $\mfg$ be the Lie algebra $\mathfrak{sl}(6,\bR)$. Recall that the notation
for the parabolic Lie subalgebra $\mathfrak p\subset \mathfrak g$ in terms
of the Satake diagrams is
%
\dyn{}{}{}{}{}, while the $\mathfrak p$-dominant highest weights as the
linear combination of the fundamental weights are often written as 
\dyn{\mu_1}{\mu_2}{\mu_3}{\mu_4}{\mu_5}.
In our preferred notation for a highest weight $\al$ this means
$\al_1=5+\mu_1+\mu_2+\mu_3+\mu_4+\mu_5$, $\al_2=4+\mu_2+\mu_3+\mu_4+\mu_5$, 
$\al_3= 3+\mu_3+\mu_4+\mu_5$, $\al_4=2+\mu_4+\mu_5$, $\al_5=1+\mu_5$,
$\al_6=0$.


We shall first discuss
two cases when the infinitesimal characters of the Verma modules are
singular. Then we come to the regular infinitesimal characters which is the most
complicated case. 

The singular infinitesimal characters can be nicely treated using
equivalence of categories due to Enright-Shelton, \cite{EnS,BH}, providing a
reduction of singular cases to the regular cases for lower rank groups.

\subsection{The case of $|2|$-singular infinitesimal characters}
The results found by Enright and Shelton, \cite{BH}, tell that the category
of Verma modules with a given $|2|$-singular infinitesimal character is
equivalent to the category of Verma modules with regular infinitesimal
character for the projective one-dimensional space $Gr(1,1)$. It is well
known that on $Gr(1,1)$, there is just one homomorphism (up to a multiple)
in each regular infinitesimal character.

The $|2|$-singular infinitesimal characters are parametrized by
$\mfg$-dominant weights in the intersection of two walls of the dominant
Weyl chamber for $\mfg$.  Hence there are six different cases parametrized
by dominant weights with two repeating components.  For each such case,
there is a unique homomorphism $\Psi$ (up to a multiple).

\begin{theorem}\label{class_2sing}
The complete list of homomorphisms $\Psi:V(\al)\to V(\be)$ between Verma modules with
$|2|$-singular infinitesimal characters is given by the lines in the following table.

$$	\begin{pmatrix}
		\la+\rho&\al&\be\\
		(aabbcd)&(abd|abc)&(abc|abd)\\
		(aabccd)&(acd|abc)&(abc|acd)\\
		(aabcdd)& (acd|abd)& (abd|acd)\\ 	
		 (abbccd)& (abc|bcd)& (bcd|abc)\\ 
		  (abbcdd)&(bcd|abd)& (abd|bcd)\\
		 (abccdd)&(bcd|acd)& (acd|bcd)	
	\end{pmatrix},
$$
where $a,b,c\in\bZ;\,a>b>c>d=0$.
\end{theorem}

\begin{proof}
The first column in the matrix is clearly giving the list (of
representatives) of all possible $|2|$-singular infinitesimal characters. 
The list of all regular infinitesimal characters for the projective space
$Gr(1,1)$ is given by $\left\{(cd)\,|\,c\in\bZ;\,c>d=0\right\}$.
Identification of both lists is given by deleting both pairs of
repeated indices in the list in the theorem.  The rest is the application
of the Enright-Shelton theory (see \cite[Theorem 7.3]{BH}).  
\end{proof}

Let us notice some familiar homomorphisms (operators going in the opposite
directions) appearing with the
lowest possible weights in this scheme. The two 1st order analogues of the Dirac
operators are: 
$$\dyn01{-4}00\to\dyn00{-3}01, \qquad \dyn00{-4}10\to\dyn10{-3}00.
$$ 
The analogue of the Laplace operator between
suitable densities has got order three: 
$$\dyn00{-4}00\to\dyn00{-2}00.$$

While in the 4-dimensional conformal case, the symbol of the conformally
invariant Laplacian (the Yamabe operator) is
given by the determinant of the 2x2 matrices viewed as $\mathfrak g_{-1}$,
here we get the same. Just the determinant is now on 3x3 matrices and thus
the order of the operator is three. 

\subsection{The case of $|1|$-singular infinitesimal characters}
By Enright and Shelton,  the category of Verma modules with a given
$|1|$-singular infinitesimal character is equivalent to the category of
Verma modules for the Grassmannian $Gr(2,2)$ with regular infinitesimal
character. The classification of all homomorphisms in such case is well
known and for every fixed infinitesimal character, it coincides with the
usual pattern with 6 standard homomorphisms and one nonstandard one (the
"long arrow"), see also \cite{SlS}. 

Recall that the standard homomorphisms appear as projections
from the Borel subalgebra case, while the non-standard ones appear in a very
different way.   

There are five different cases parametrized by dominant
weights with one repeating component.  All homomorphisms in each pattern are
unique (up to a constant multiples).

\begin{theorem}\label{class_1sing}
The complete list of homomorphisms between Verma modules with $|1|$-singular
infinitesimal character is given by the following list. All five cases
are forming the same pattern:

\smallskip	 
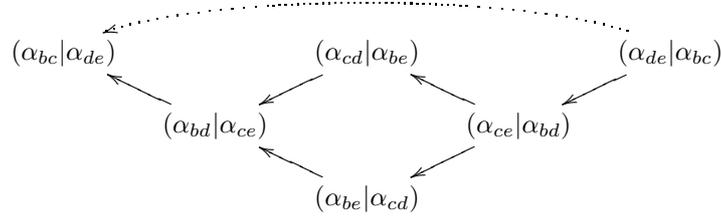
\begin{equation*}
\xymatrix@R=4mm@C=4mm{
		\wgt{\al_{bc}}{\al_{de}}
		&&\wgt{\al_{cd}}{\al_{be}}
		\ar[dl]
		&& \wgt{\al_{de}}{\al_{bc}}
		\ar[dl]\ar@(ul,ur)@{..>}[llll]
		\\
		&\wgt{\al_{bd}}{\al_{ce}}
		\ar[ul]
		&&\wgt{\al_{ce}}{\al_{bd}} 
		\ar[ul]\ar[dl]
		&\\
		&&\wgt{\al_{be}}{\al_{cd}}
		\ar[ul]&&  
}
\end{equation*}
with weights given by the lines in the table
$$	\begin{pmatrix}
	\la+\rho&(\al_{bc}|&(\al_{bd}|&(\al_{be}|&(\al_{cd}|&(\al_{ce}|&(\al_{de}|\\
	(aabcde)&(abc|ade)&(abd|ace)&(abe|acd)&(acd|abe)&(ace|abd)&(ade|abc)\\
(abbcde)&(abc|bde)&(abd|bce)&(abe|bcd)&(bcd|abe)&(bce|abd)&(bde|abc)\\
(abccde)&(abc|cde)&(acd|bce)&(ace|bcd)&(bcd|ace)&(bce|acd)&(cde|abc)\\ 	
(abcdde)&(abd|cde)&(acd|bde)&(ade|bcd)&(bcd|ade)&(bde|acd)&(cde|abd)\\
(abcdee)&(abe|cde)&(ace|bde)&(ade|bce)&(bce|ade)&(bde|ace)&(cde|abe)	
\end{pmatrix},
$$
where $a,b,c,d\in\bZ;\,a>b>c>d>e=0.$

Full arrows in the diagram are the standard homomorphisms, 
which form the (generalized) BGG resolution.
The dotted arrow is a nonstandard homomorphism. The only nontrivial
compositions appear in the central diamond and they are equal up to sign. 
\end{theorem}

\begin{proof}
The first column in the matrix is clearly giving the list (of
 representatives) of all possible $|1|$-singular infinitesimal characters. 
 The list of all regular infinitesimal characters for the Grassmannian
 $Gr(2,2)$ and homomorphisms between them is well known and it is given in the
 theorem. The identification of both lists and patterns is given by the deleting 
the pair of repeating
indices.  The rest is again the application of
 the Enright-Shelton theory (see \cite[Theorem 7.3]{BH}).
\end{proof}
 
Again, notice the analogue of the invariant square of the Laplacian (the Paneitz
operator) between suitable densities here. It is the nonstandard operator in 
the middle line
with the lowest possible weights and its order is six, while its symbol is
(again in full analogy to the 4-dimensional case) the square of the
determinant: 
$$
\dyn00{-5}00\to\dyn00{-1}00
.$$

\subsection{The case of regular infinitesimal character}
The list of all homomorphisms between Verma modules with regular infinitesimal
characters
can be read off \cite{BC} for low rank semisimple groups.\footnote{It also 
can be effectively computed by Brian Boe's
software, whose Fortran77 implementation the authors gratefully 
use for more than two decades already.} 

\begin{theorem}\label{class_regular}
The list of all homomorphisms between Verma modules with 
a fixed regular infinitesimal character $(abcdef)$
is given in the following diagram,
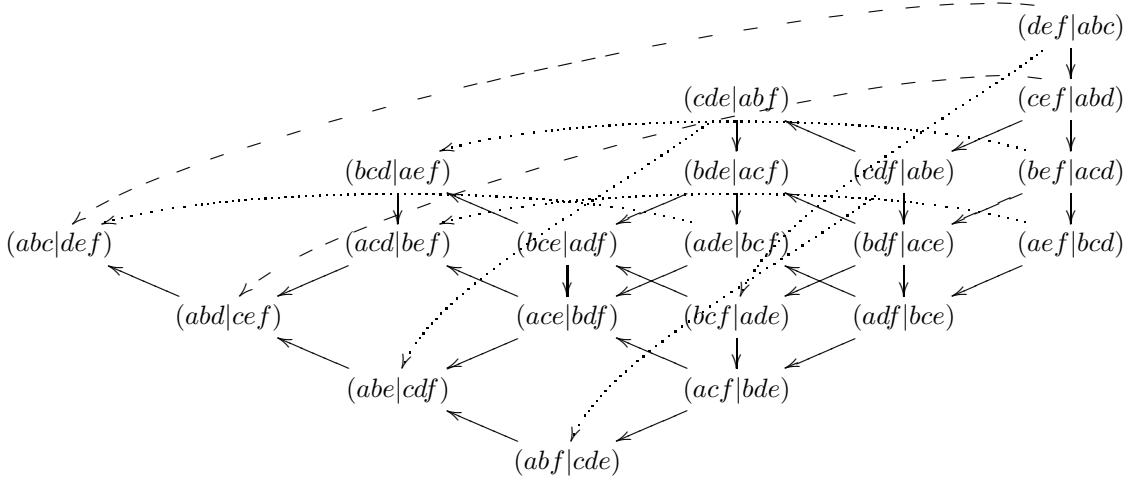
\begin{equation*}
\xymatrix@R=4mm@C=6mm{
&&&&&& \wgt{def}{abc} 
\ar@(u,u)@{-->}[llllllddd]
\ar[d]
\ar@(dl,u)@{..>}[ddddll]
\\
&&&& \wgt{cde}{abf} 
\ar[d] 
\ar@(dl,u)@{..>}[ddddll]
&& \wgt{cef}{abd}
\ar[d] \ar@[green][dl]
\ar@(u,u)@{-->}[lllllddd]
\\
&& \wgt{bcd}{aef}
\ar[d]
&& \wgt{bde}{acf}
\ar@[blue][ld] \ar[d] 
& \wgt{cdf}{abe}
\ar@[green][ul] \ar[d] \ar@(dl,u)@{..>}[ddddll]
&\wgt{bef}{acd}
\ar@(ul,ur)@{..>}[llll]
\ar@[blue][dl] \ar[d]
\\
\wgt{abc}{def}
&& \wgt{acd}{bef}
\ar@[red][dl] 
& \wgt{bce}{adf}
\ar@[blue][ul] \ar[d] 
& \wgt{ade}{bcf}
\ar@[red][dl] \ar@(ul,ur)@{..>}[llll]
& \wgt{bdf}{ace}
\ar@[blue][ul] \ar[d] 
\ar@[blue][dl]  
& \wgt{aef}{bcd}
\ar@[red][dl] \ar@(ul,ur)@{..>}[llll]
\\
& \wgt{abd}{cef}
\ar@[red][ul] 
&& \wgt{ace}{bdf}
\ar@[red][ul] \ar@[red][dl] 
& \wgt{bcf}{ade}
\ar[d] \ar@[blue][ul]
& \wgt{adf}{bce}
\ar@[red][ul] \ar@[red][dl]
\\
&& \wgt{abe}{cdf}
\ar@[red][ul] 
&& \wgt{acf}{bde}
\ar@[red][ul] \ar@[red][dl]
\\
&&& \wgt{abf}{cde}
\ar@[red][ul] 
} 
\end{equation*}
%
%
where $a,b,c,d,e\in\bZ;\ a>b>c>d>e>f=0$.
Full arrows denote standard homomorphisms (included in the corresponding BGG resolution). 
Dotted and dashed arrows 
describe nonstandard homomorphisms.
 
For every square of arrows, the homomorphism given by the diagonal composition is also nontrivial.
This includes the cases of squares given by standard homomorphisms 
as well as those given by two standard and two nonstandard homomorphism. On
top of that, the diagonal in the cube in the middle of the diagram (a
composition of three morphisms) is also nontrivial.  

In the lowest regular infinitesimal character, the standard homomorphisms correspond
to the relevant restrictions and projections of the exterior differential on
exterior forms, further there are 6 nonstandard operators of
order 4 (dotted arrows), one nonstandard operator of order 7, and one
nonstandard operator between functions and top forms (of order 9).  The
structure of the scheme is independent of the chosen infinitesimal
character, but the order of operators is changing.  
\end{theorem}

\begin{proof}
Regular infinitesimal characters of Verma modules for $Gr(3,3)$ are parametrized
equivalence classes of weights of the form $(abc|def);\;a,b,c,d,e,f\in\bZ;\;a>b>c>d>e>f$
modulo the weight $(111111).$ 
The action of the Weyl group is given by permutation of their components but we keep just
those dominant for $\mfp.$
Representatives of equivalence classes are chosen
to be weights with $f=0$.

The lowest infinitesimal character is given by $(543210)$. 
The corresponding nodes in the scheme
are given by irreducible pieces of the Grassmann 
algebra $\Lambda^*(\mfg_1)$ with respect to the action of the Levi 
subalgebra $\mfg_0$.
 
The list of all nonstandard homomorphisms for the lowest infinitesimal character can be found
in \cite[p. 82]{BC}. 
It is possible to check there all cases of nontrivial compositions 
of operators in squares and the one in the cube.

The scheme is repeated for more complicated infinitesimal characters and corresponding homomorphisms
are obtained by the standard translation procedure, see the comments on the
Jantzen-Zuckermann translation principle below.
The shifts of the scheme to other infinitesimal characters have the same structure.
\end{proof}


The big cell in the Grassmannian is isomorphic to the space of real matrices
$(x_{ij})_{i,j=1}^3$ and symbols of homomorphisms corresponding to
intertwining operators on line bundles can be directly understood as 
polynomials in terms of the entries of these matrices. 
In particular, there is the series of intertwining differential
operators on line bundles with symbols given by $\left(\det(\frac{\pa}{\pa
x_{ij}})\right)^n$, for all $n\in\bN$. They are the analogues of powers
of Laplacians in conformal geometry (see \cite{HS}).  The corresponding
 homomorphisms are easily visible in our classification. We have seen already that the $n=1$ case appears in the
infinitesimal character $(322110)$, while the $n=2$ case is in the infinitesimal character $(432210)$. For $n\geq 3$, they are the longest dashed arrows in the 
regular infinitesimal characters $(n+2,n+1,n,2,1,0)$. In the other notation,
they appear as
$$
\dynkin\root{0}\link\root{0}\link\nodroot{-n-3}\link\root{0}\link\root{0}\enddynkin
\to\dynkin\root{0}\link\root{0}\link\nodroot{n-3}\link\root{0}\link\root{0}\enddynkin
.$$

\section{The translation principle}

Now we are coming to the crucial point for understanding the overall
structure of the homomorphisms, a functorial construction providing
(iso)morphisms between the categories of modules with given
infinitesimal characters. 

We refer to the survey in \cite{SlS} and the references therein, here we
recall the construction very briefly. The main observation is very simple:

Suppose $\mathbb W$ is a $G$-module (and $P$-module by restriction), and
$\mathbb E$ a $P$-module. Then 
$\mathfrak U(g)\otimes \mathbb E^*\otimes \mathbb W^*$, as a
$\mathfrak g$-module, can be viewed in two
different ways:
\begin{gather*}
X(x\otimes e \otimes w) = Xx\otimes e\otimes w
\\
X(x\otimes e \otimes w) = Xx\otimes e\otimes w + x\otimes e \otimes Xw
,\end{gather*}
for all $X\in\mathfrak g$, $x\in\mathfrak U(\mathfrak g)$, $e\in\mathbb E^*$,
and $w\in\mathbb W^*$.
While the first case descends to the $(\mathfrak
U(\mathfrak g),P)$-module structure on $V(\mathbb
E\otimes\mathbb W)$, the second one leads to the structure of $V(\mathbb
E)\otimes \mathbb W^*$. The identity on $1\otimes e \otimes w$
 extends to an isomorphism between these two modules:
\begin{equation}\label{twisting}
V(\mathbb E\otimes \mathbb W) \simeq V(\mathbb E)\otimes \mathbb W^*
.\end{equation}

Non-trivial irreducible $G$-modules 
$\mathbb W$ are never
irreducible $P$-modules. They always enjoy a composition series 
\begin{equation}\label{composition_series}
\mathbb W = \mathbb W_{\al-\ell} + \mathbb W_{\al-\ell+1} + \dots + \mathbb
W_{\al}
\end{equation} 
where the labeling reflects the scalar action by the grading element in
$\mathfrak g_0$, and
the `right ends' $W_j=\mathbb W_j+\dots+\mathbb W_\al$ form $\mathfrak p$-submodules, i.e., we get the filtration
\begin{equation}\label{filtration}
\mathbb W_\al = W_\al\subset W_{\al-1}\subset\dots\subset
W_{\al-\ell}=\mathbb W,
\end{equation}
with $\mathbb W_j=W_j/W_{j+1}$. 
As a $\mathfrak g_0$-module,
the composition series is a direct sum of submodules $\mathbb W_j$ and 
each of them further decomposes into $\mathfrak g_0$-irreducible submodules 
$\mathbb W_{j,k}$. 

For each irreducible $\mathfrak p$-module $\mathbb E$, and a
$\mathfrak g$-module $\mathbb
W$ as above,
$$
\mathbb E\otimes \mathbb W = \mathbb E\otimes \mathbb W_{\al-\ell}+ \dots +
\mathbb E\otimes \mathbb W_\al
,$$
and all components $\mathbb E\otimes \mathbb W_i$ split into direct sum of
irreducible $\mathfrak g_0$-modules $\mathbb E_{i,j}$.  

If we assume that one of the many modules $\mathbb E'=\mathbb E_{i,j}$
has got a distinct infinitesimal character than all the other modules 
in the above decomposition, then the injection $V(\mathbb E')\to
V(\mathbb E\otimes \mathbb W)$ is defined by its image being the joint
eigenspace of the infinitesimal character of $V(\mathbb E')$. Clearly,
there is the complementary subspace defined as the generalized eigenspaces
of all the other infinitesimal characters there. 

Thus, under the latter assumption, $V(\mathbb E')$ canonically splits
off the $V(\mathbb E\otimes \mathbb W)=V(\mathbb E)\otimes \mathbb W^*$ as a
direct summand, and this provides the translation idea: 

For each non-trivial $(\mathfrak U(\mathfrak g),P)$-module homomorphism 
$\Phi:V(\mathbb F)\to V(\mathbb
E)$, and assuming further that both $V(\mathbb E')$ and $V(\mathbb
F')$ enjoy the same and unique infinitesimal character in $V(\mathbb E\otimes
\mathbb W)$ and $V(\mathbb F\otimes \mathbb W)$, respectively, we
obtain the composed homomorphism
\begin{equation}\label{translation}
V(\mathbb F')\to V(\mathbb F\otimes \mathbb W)=V(\mathbb F)\otimes \mathbb W^*\to
V(\mathbb E)\otimes \mathbb W^*= V(\mathbb E\otimes \mathbb W) \to V(\mathbb
E')
.\end{equation}
We talk about \emph{twisting the homomorphism} $\Phi$ by tensoring it with the identity
on $\mathfrak g$-module $\mathbb W^*$. Although it is highly non-trivial to 
recognize whether the \emph{translated 
morphism} $V(\mathbb F')\to V(\mathbb E')$ is nontrivial, there are general
arrangements securing this. In particular, the celebrated 
Jantzen-Zuckermann translation principle tells, that using the highest and
lowest weights of $\mathbb W$ always provides two adjoint equivalences of
categories for all equisingular infinitesimal characters. 

Our aim is to extend this algebraic translations to the realm of
curved Cartan geometries, but not all $\mathfrak g$-modules are generally
suitable for this. Thus, the crucial
observation in \cite{ES} was that many other weights in $\mathbb W$ work as
well. We provide two such observations as two propositions, both explained
carefully in \cite{SlS}:

\begin{proposition}[\cite{ES, SlS}]\label{transprop1}
Suppose that $V({\Bbb E})$ and $V({\Bbb F})$ have the same
infinitesimal character. Suppose that $V({\Bbb E}^\prime)$ and $V({\Bbb
F}^\prime)$
have the same infinitesimal character. Let ${\mathbb W}$ be a 
finite-dimensional
irreducible representation of $G$ and suppose that
\begin{itemize}
\item $V({\Bbb F}^\prime)$ occurs in the composition series for
      $V({\Bbb F}\otimes{\mathbb W})$ and has distinct infinitesimal 
character from all other factors;
\item $V({\Bbb E}^\prime)$ occurs in the composition series for
      $V({\Bbb E}\otimes{\mathbb W})$ and has distinct infinitesimal 
character from all other factors.
\end{itemize}
It follows that $V({\Bbb F})$ occurs in the composition series for
$V({\Bbb F}^\prime\otimes{\mathbb W}^*)$ and that $V({\Bbb E})$ occurs in the
composition series for $V({\Bbb E}^\prime\otimes{\mathbb W}^*)$. We suppose
further that
\begin{itemize}
\item all other composition factors of $V({\Bbb F}^\prime\otimes{\mathbb
W}^*)$ have infinitesimal character distinct from $V({\Bbb F})$;
\item all other composition factors of $V({\Bbb E}^\prime\otimes{\mathbb
W}^*)$ have infinitesimal character distinct from $V({\Bbb E})$.
\end{itemize}
Then translation gives an isomorphism
$$\operatorname{Hom}_{({\frak U}({\mathfrak g}),P)}(V({\Bbb F}),V({\Bbb
E}))\simeq
  \operatorname{Hom}_{({\frak U}({\frak g}),P)}(V({\Bbb F}^\prime),V({\Bbb E}^\prime))$$
whose inverse is given by translation using~${\mathbb W}^*$.
\end{proposition}

The other proposition deals with the one-way translations producing
non-trivial homomorphisms in the less singular patterns from the more
singular ones. 
Recall we write $\varphi(\mathbb E)$ for the constant action of the grading
element.

\begin{proposition}[\cite{ES,SlS}]\label{transprop2}
Let $\Phi:V(\mathbb E)\to V(\mathbb F)$ be a nontrivial homomorphism of 
Verma modules,
and let $\mathbb W$ be an irreducible finite dimensional  $G$-module.

Suppose  that there are irreducible $\mathfrak p$-modules
$\mathbb E_1,\mathbb E_2,\mathbb F_1,\mathbb F_2$ such that:
	
\noindent (i) 
	$\mathbb E\otimes\mathbb W=\mathbb E_1\oplus
\mathbb E_2\oplus\mathbb E';\;\mathbb F\otimes\mathbb W=\mathbb F_1\oplus \mathbb F_2\oplus \mathbb F';$
	
\noindent (ii) Verma modules $V(\mathbb E_1),V(\mathbb E_2),V(\mathbb F_1), V(\mathbb F_2)$ have the
same  infinitesimal character;
	
\noindent (iii)
	all pieces in the composition series for $V(\mathbb E'), V(\mathbb F')$ have
different infinitesimal characters.
	from those in (ii);
	
\noindent (iv) $ \varphi(\mathbb E_1)< \varphi(\mathbb E_2)$,
$\varphi(\mathbb F_1)>\varphi(\mathbb F_2)$  
 and  
	$V(\mathbb F)$ splits off from   $V(\mathbb F_1\otimes\mathbb W^*).$ 
	
	If there is no nontrivial homomorphism from $V(\mathbb E_2)$  to
$V(\mathbb F_1)$, 
	then the translated homomorphism
	$$\hat{\Phi}: V(\mathbb E_1)\to V(\mathbb E\otimes\mathbb W)\to  V(\mathbb F\otimes\mathbb W)\to
V(\mathbb F_1)$$ is nontrivial.    
\end{proposition}

%
%

Let us come back to our $Gr(3,3)$ situation. The next lemmas show,
that we have got nice and straightforward possibilities to twist
homomorphisms with simple modules. Notice, we are writing the $\mathfrak
g$-dominant weight $\la$ for the $\mathfrak g$ modules, while keeping the
$\al$ notation for the $\mathfrak p$-dominant ones (i.e., shifted by the
lowest weight $\rho$).

\begin{lemma}\label{lem_100000}
Consider a $\mathfrak p$-dominant weight $\al=(abc|def)$, 
$a>b>c,\,d>e>f$.

\noindent (i) Its tensor product with $\mathbb W=(100000)$ decomposes as
 	\begin{equation}\label{dec}
 		\begin{split}
 			(abc|def)\otimes (100000)=(a+1,b,c|def)+
(a,b+1,c|def)+(a,b,c+1|def) + {} 
\\ (abc|d+1,e,f)+(abc|d,e+1,f)+(abc|d,e,f+1)	
 		\end{split}
 	\end{equation}	
supposing that summands which are not $\mathfrak p$-dominant 
are removed from the right-hand side.
 	
\noindent (ii) If $\al$ has got regular infinitesimal character, then all 
summands in (\ref{dec}) have pairwise different 
infinitesimal characters.

\noindent (iii) If the infinitesimal character of $\al$ is singular, then the
terms on the right hand side increasing the non-repeating values have got an
infinitesimal character different from the rest. 

\noindent (iv) Tensor products with $\mathbb W^*=(111110)$ behave similarly
and can be understood as decreasing one of the values in $\al$. 
\end{lemma}

\begin{proof}
(i) The relation (\ref{dec}) follows from
$ 	(100000)=(100|000)	+ (000|100)	$
 	and from $\mathfrak{sl}(3,\mathbb R)$ decomposition
$$(abc)\otimes (100)=(a+1,b,c)\oplus (a,b+1,c)\oplus (a,b,c+1),
$$
where $a>b>c$ and the summands are not there if components are not pairwise different.
 	
(ii) The argument is seen well on an example:   
 $(a+1,b,c|def)$ and $(a,b,c|d+1,e,f)$ have four components in common 
and they have
the same infinitesimal character if and only if the pair $(a+1,d)$ equals 
to the pair $(a,d+1)$.
But this is true only if $a=d$ which is excluded in regular character.
The same reasoning works for any pair in the sum.

(iii) The above reasoning fails if $a=d$, but we still conclude the same
claim if $a$ does not belong to the repeated values.

(iv) In this case we run through the same procedure, just going through the
options for choosing one of the values which is not increased by one. By the
equivalence relation on the weight vectors modulo (111111) this corresponds
to subtracting one from the chosen value.
\end{proof}

\begin{lemma}\label{lem_110000}
Consider a weight $\al=(abc|def)$, $a>b>c$, $d>e>f$, with a singular
infinitesimal character. Then its tensor product with $\mathbb W=(110000)$
contains summands with one of the couples of repeated indices increased both by
one (whenever this increase makes sense) and 
this summand has got distinct infinitesimal character from the rest. 

The tensor product with $\mathbb W*=(111100)$ works similarly, and it can be
understood as decreasing the couples repeated values in the weights by one,
resulting in infinitesimal characters distinct from the rest. 
\end{lemma}

\begin{proof} Here $\mathbb W$ has got the composition series 
$(110|000)+(100|100)+(000|110)$. Its tensor product with $\al$ has got 15
potential
summands, but only those coming from the middle part are of interest for us
now. Clearly they involve all possible increases of the repeated couples by
one. A similar direct check as above reveals that the resulting weights have got distinct
infinitesimal character from the rest.

The claim about $\mathbb W^*$ is deduced the same way.
\end{proof}

Next, we show how to (uniformly) translate individual components of the
infinitesimal character $(abcdef)$ by one.

\begin{proposition}\label{simple_trans}
Twisting the morphisms with $\mathbb W=(100000)$ and $\mathbb W=(110000)$
allows to translate the entire schemes of equisingular weights bijectively.

In particular, we can obtain all the diagrams of morphisms from Theorems
\ref{class_2sing}, \ref{class_1sing}, \ref{class_regular} from the diagrams
with the lowest possible weights by these translations. 
\end{proposition}

\begin{proof}
The general reasoning is again best seen on examples. 
Suppose that $b>c+1$ and that we want to move by translation
the whole Hasse diagram from regular character $(abcdef)$ 
to character $(a,b,c+1,d,e,f)$. 
We can do that by tensoring every term in the Hasse diagram by 
$(100000)=(100|000)+(000|100)$
and restrict everything to the character $(a,b,c+1,d,e,f)$.

Indeed, we know from Lemma \ref{lem_100000} 
that all summands (at any position in the Hasse diagram) 
have pairwise different  characters. 
But at any position in the Hasse diagram, there clearly is
a unique summand sitting in the character $(a,b,c+1,d,e,f)$.
At the same time, the tensor products of the results with $\mathbb W^*$
behave well, again by Lemma \ref{lem_100000}, so we may exploit Proposition
\ref{transprop1}. Summarizing, 
by restricting everything to character of $(a,b,c+1,d,e,f)$,
we have performed the (uniform) translation of the whole diagram.

Dealing with singular infinitesimal characters we either use the same
procedure as above (if we want to increase one of the non-repeating values),
or we have to be more careful and use Lemma \ref{lem_110000}, which
again leads to applications of Proposition \ref{transprop1}.
\end{proof}


Thus, we have verified that all equisingular
patterns for our Grassmannian geometries always have identical shapes.
Moreover, notice that we are always using a $\mathfrak g$-module with the
composition series of length at most three, and so the injection and projection
operators in the translations are always of order at most two!

\section{Semiholonomic Verma modules and morphisms}

While the invariant differential operators between homogenous bundles over
Klein geometries $G/P$ can be completely understood via the purely algebraic
and functorial setup mentioned above, the curved Cartan geometries are not
that easy. 

The homogeneous bundles determined by $P$-modules extend to functors of the
associated bundles to the Cartan bundles $\mathcal G$, but the
question which of the invariant operators extend to natural differential
operators between such bundles is much more intriguing.  On one hand, the
standard operators in the schemes with regular infinitesimal characters
always admit such extensions and they are available by direct constructions,
known as the BGG machinery, cf.  \cite{CSS, CD}.  At the same time, in
conformal Riemannian geometry, there is a class of operators which do not
allow such extensions at all, see \cite{ES, SlS} for explanation.

The approach introduced in \cite{ES} allows to mimic the Klein geometry
algebraic story and deliver a lot of the non-standard operators in a purely
algebraic way. 
We shall follow this track for our almost Grassmannian geometries modelled
over $Gr(3,3)$ now.

The main tool consists in understanding the semi-holonomic jets of sections
of the associated bundles $\mathcal E$
corresponding to a $P$-module
$\mathbb E$ and the universal differential operator $\sigma \mapsto \bar
j^k\sigma$ on sections $\sigma$ of the bundles $\mathcal E$, 
provided by the iterated Cartan connection fundamental
derivative. The analogous construction to the Verma modules then provides
the so-called semi-holonomic universal enveloping algebra $\bar {\mathfrak
U}(\mathfrak g)$ and the topological dual to the infinite semi-holonomic jets
$\bar V(\mathbb E) = \bar{\mathfrak U}\otimes_{\mathfrak U(\mathfrak
p)}\mathbb E^*$.

The latter $\mathfrak g$-modules are again highest weight modules and 
there is the following commutative 
diagram of short exact sequences:
\begin{equation}\label{jetdiagram}
\xymatrix@R=5mm@C=10mm{
0 \ar[r] 
& \bar V_{k-1}(\mathbb E) \ar[r] \ar[d]
& \bar V_k(\mathbb E) \ar[r] \ar[d] 
&\otimes^k(\mathfrak g/\mathfrak h)\otimes \mathbb E^* \ar[r] \ar[d]
& 0\\
0 \ar[r] 
& V_{k-1}(\mathbb E) \ar[r]
& V_k(\mathbb E) \ar[r]  
& S^k(\mathfrak g/\mathfrak h)\otimes \mathbb E^* \ar[r]
& 0
}
\end{equation}
where the most right vertical arrow is given by the symmetrization.
 
All this is carefully explained in \cite{SlS} and we are getting the
opportunity to study the possible liftings of
the existing homomorphisms $V(\mathbb F)\to V(\mathbb E)$ to morphisms 
$\bar V(\mathbb F)\to \bar V(\mathbb E)$ with respect to the canonical
projection. Moreover, due to the Frobenius reciprocity, this is equivalent
to the search for the dashed $P$-module morphisms in the following commutative diagram:
\begin{equation}\label{liftingdiagram}
\xymatrix@R=5mm@C=20mm{
& \bar V(\mathbb E) \ar[d]
\\
\mathbb F^* \ar@{-->}[ur] \ar[r] 
& V(\mathbb E) 
}
\end{equation}
In turn, for irreducible modules $\mathbb E$, $\mathbb F$, 
this is equivalent to finding a highest weight vector in $\bar
V(\mathbb E)$ covering the relevant highest weight vector in $V(\mathbb
E)$.
The order of a
homomorphism $\Phi$ between the semi-holonomic Verma modules is defined in the same
way as for the standard Verma modules, and the
symbol of $\Phi$ is
$$
\sigma(\Phi):\mathbb F^* \to \bar V_k(\mathbb E)\to \bar V_k(\mathbb E)/\bar
V_{k-1}(\mathbb E) = \otimes^k(\mathfrak g_-)\otimes\mathbb E^*
.$$ 

Due to the existence of the universal differential operators $\bar
j^k:\mathcal E \to \bar J^k\mathcal E$, every nontrivial homomorphism $\Phi:\bar
V(\mathbb F)\to \bar V(\mathbb E)$ obviously provides an invariant
differential operator $\mathcal E \to \mathcal F$, 
of the same order.  

The final ingredient we need is the following simple fact:

\begin{proposition}[\cite{ES,SlS}]
A homomorphism $V(\mathbb F)\to V(\mathbb E)$ of order at most two 
always lifts to a homomorphism $\bar V(\mathbb F)\to \bar V(\mathbb E)$.
\end{proposition}

This simple proposition implies again that all first and
second order linear invariant differential operators extend canonically from the
Klein's models to the entire category of the corresponding Cartan
geometries. This concerns also the initial embeddings and final projection
in the translation construction, whose order is equal (for all $|1|$-graded
geometries) to the distance of the components from the highest or lowest
weight components in the tensor product composition series, respectively.

Finally, 
we can recover the translation idea from the Klein's model in the curved case:
$$
\xymatrix@R=5mm@C=6mm{
\bar V(\mathbb F') \ar[r] \ar[d] \ar@(ur,ul)@{-->}[rrrrr]^{\bar \Phi'}
&\bar V(\mathbb F\otimes \mathbb W) \ar@{=}[r] \ar[d]
&\bar V(\mathbb F)\otimes\mathbb W^* \ar[d] \ar[r]^{\bar \Phi\otimes 1}
&\bar V(\mathbb E)\otimes \mathbb W^* \ar@{=}[r] \ar[d]
&\bar V(\mathbb E\otimes \mathbb W) \ar[r] \ar[d]
&\bar V(\mathbb E') \ar[d]
\\
V(\mathbb F') \ar[r] \ar@(dr,dl)@{-->}[rrrrr]_{\Phi'}
&V(\mathbb F\otimes \mathbb W) \ar[r]
&V(\mathbb F)\otimes\mathbb W^* \ar[r]^{\Phi\otimes 1}
&V(\mathbb E)\otimes \mathbb W^* \ar[r] 
&V(\mathbb E\otimes \mathbb W) \ar[r] 
&V(\mathbb E')
}
$$

Consequently, starting with any homomorphism $\Phi$ of Verma modules which lifts
to $\bar \Phi$,
all other morphisms $\Phi'$ obtained from $\Phi$ 
by means of translations described in Propositions \ref{transprop1}
and \ref{transprop2}, i.e., using only splitting operators of order at most two,
admit the covering by homomorphisms $\bar \Phi'$ of semi-holonomic Verma
modules, too.
Moreover, the symbols of the covered homomorphisms are obtained by
symmetrizations of symbols of those covering ones.  
Again, see \cite{SlS} for detailed explanation.


We are ready to formulate the main results of this paper now. Both theorems
deal with the case of the almost Grassmannian geometries modelled over
$Gr(3,3)$, i.e., $\mathfrak g=\dyn{}{}{}{}{}{}$. 

\begin{theorem} \label{th6}
All homomorphisms between Verma modules with singular infinitesimal characters lift
to homomorphisms between the corresponding semi-holonomic Verma modules.
\end{theorem} 

\begin{theorem}\label{th7}
All homomorphisms between Verma modules with regular 
infinitesimal characters $(abcdef)$, $a>b>c>d>e>f$, lift
to homomorphisms between the corresponding semi-holonomic Verma modules,
except the following three classes:
$$
V(def|abc)\to V(abc|def);
V(ade|bcf)\to V(abc|def);
V(def|abc)\to V(bcf|ade).$$
These three classes of homomorphisms do not allow any lift.
\end{theorem} 

In the theorem, there are three sets of homomorphisms, which
are excluded from the statement. We prove that the second and third cases do
not allow for any lifts in Theorem \ref{no_lift} below, while the first case, the longest arrow, is resolved by quite different methods in the last section.
 
 
In order to see that all homomorphisms in the same position in our scheme of
equisingular infinitesimal characters do not allow for lifts, it
is enough to deal with one of them. Indeed, the general simple translation
procedure for equisingular infinitesimal characters from Proposition
\ref{simple_trans} ensures that the entire
class must behave equally. 

The next theorem resolves the two nonstandard homomorphisms of order $4$ in 
the lowest possible infinitesimal character and this proves a part of the very last
claim in Theorem \ref{th7}. 
 
\begin{theorem}\label{no_lift}
The homomorphisms 
\newline 	
 	(A) $(521|430)\to (543|210)$	
\newline	
 	(B) 	 $(210|543)\to (430|521)$
\newline	
cannot be lifted to semi-holonomic case.
\end{theorem}
 
 \begin{proof}
Note that there are $9$ roots in $\mfg_{-1}$. Let us denote by 
$y_{ij}$, $3\leq i\leq 5$, $0\leq j\leq 2$,
the matrices having $1$ at the position $(ij)$ and zero everywhere else.
They are root vectors for the roots $\varepsilon_i-\varepsilon_j$ and they 
form a basis for $\mfg_{-1}$.

(A)
In terms of highest weights, we deal with the homomorphism
$V(\bE)\to V(\bF)$,
 	where $\bE=(0,-2,-2|2,2,0)$ and $\bF= (000|000).$ Let us consider
 	Lie subalgebra of $\mfg$ given by matrices in $\mfg$ having the
 	first and the last rows, as well as the first and 
last columns trivial.  This
 	is a Lie algebra $\mfg'$ of type $A_3$.  The intersection $\mfp'$ of
 	the parabolic subalgebra $\mfp$ of $\mfg$, corresponding to the
 	middle node crossed with $\mfg'$, is clearly a parabolic subalgebra
 	of $\mfg'$ and the corresponding geometry is the 4-dimensional conformal
 	geometry.
 	
 	Suppose now that there is a lift of the homomorphism $(A)$ to
 	semi-holonomic Verma modules corresponding to a singular vector
 	$w\otimes e$, where $e$ is a basis vector for the trivial module
 	$\bF$, and $w$ has the weight $(0,-2,-2|2,2,0)$, and it is a
 	(noncommutative) polynomial in variables $y_{ij}\in \mfg_{-1}$. Due
 	to its weight, the polynomial depends only on $y_{ij}$, $3\leq
 	i\leq 4$, $1\leq j\leq 2$. 
Then, the same formula also defines a singular
 	vector for a semi-holonomic lift of a nonstandard homomorphisms
 	in $\mfg'$, because the actions of elements in $\mfg'_{0}$ and
 	$\mfg'_1$ on $w\otimes e$ coincide with the actions of elements in
 	$\mfg_{0}$ and $\mfg_1$ on $w\otimes e$. This is a contradiction
 	with the well known results in 4-dimensional conformal geometry
 	(cf. \cite{ES,SlS}).
 	
 	(B)
Now we deal with the homomorphism
 	$V(\bE)\to V(\bF)$,
 	and weights $\bE=(-3,-3,-3|3,3,3)$ and $\bF= (-1,-1,-3|3,1,1)$. 
The same reasoning also works here.
 	
This time, the subalgebra $\mfg'$ is given by matrices from $\mfg$ 
having trivial
 	the third and the fourth rows and the third an the fourth columns.
 	Again, the resulting Lie algebra is of type $A_3$.
 	
 	Suppose now that there is a lift of the homomorphism (B) to semi-holonomic Verma modules
 	corresponding to a singular vector in $V(\bF)$.
 	The representation $\bF$ is no more one-dimensional, 
hence the expression for the
 	singular vector is more complicated one but due to the weight of the
 	singular vector, only allowed action of lowering operators in the
 	formula is by elements from $\mfg'_{-1}$ and so the formula
 	recovers a singular vector in 4-dimensional conformal geometry.
 \end{proof}
 
\section{The proof of Theorem  \ref{th6}}

Let us recall the two Propositions \ref{transprop1}, and \ref{transprop2}. In
particular, the latter one allows us to employ suitable one-way translations
from more singular to less singular patterns. Dealing with the semi-holonomic
Verma modules, we prefer to translate with $\mathfrak g$-modules
with composition series of length at most three since then all the
projections and injections in the twisting procedure lift.

Concerning the notation for weights, we shall again use the dominant $\mathfrak
g$-weights $\la$ (i.e., without adding the lowest weight $\rho$), while we
keep writing $\al$ for the $\mathfrak p$-dominant weights. This might look
chaotic, but it is easier to see the tensor product decompositions this way.

 \subsection{The case of $|2|$-singular Verma modules}
 (i) Let us first consider the lowest $|2|$-singular infinitesimal characters. All homomorphisms
 are of first or second order (hence they lift) with the exception of
the homomorphism  $(210|321)\to(321|210)$,
which is of the third order. 

Using Proposition \ref{transprop1}, we shall get it as a translation of the 
homomorphism  $(210|310)\to(310|210)$ by $\mathbb W=(110000)$ (the most relevant
weights are boxed):
\begin{align*}
\mathbb F &= (310|210);\ \mathbb E = (210|310);\ \mathbb F' = (321|210);\
\mathbb E' = (210|321)
\\
\mathbb F\otimes \mathbb W &= (420|210) + \boxed{(321|210)} + (410|310) +
(320|310)  + (310|320)
\\
\mathbb E\otimes\mathbb W &= (320|310) + (310|410) + (310|320) + (210|420) +
\boxed{(210|321)}
\\
\mathbb F'\otimes \mathbb W^* &= (432|310) + (431|320) + \boxed{(421|321)}
\\
\mathbb E'\otimes\mathbb W^* &= \boxed{(321|421)} + (320|431) + (310|432)
\end{align*}
and the conditions from the proposition are easy to check.

As we learned in Proposition \ref{simple_trans}, 
by further translations using Proposition \ref{transprop1}, 
with $\mathbb W=(100000)$ or $W=(110000)$, we can then reach
any homomorphism in all equisingular infinitesimal characters.
 
\subsection{The case of $|1|$-singular Verma modules}
We start first with lowest
weight involving components $(43210)$ with one of them repeated.
All standard homomorphisms (appearing as arrows in the BGG complex for the
four-dimensional conformal Riemannian case) 
are in this case of order at most two,
hence they lift.

We are going to obtain one of the remaining nonstandard operator by a one-way 
translation (see Proposition \ref{transprop2}) from a $|2|$-singular case:
$$ V(210|321)\to V(321|210)\implies V(210|431) \to V(431|210)).$$

Notice, the third order morphism on the left is the 3rd order analog of
Laplace which lifts, as we already proved
above. The target of the translation is the 5th order morphism in the 4th
line in Theorem \ref{class_1sing}.

Consider modules
\begin{align*}
	\bF&=(321|210);\ 	\bF_1=(431|210);\ 	\bF_2=(421|310)\\
	\bE&=(210|321) ;\  \bE_1=(210|431);\  \bE_2=(310|421)\\
	\bW&=(110000)=(110|000)+(100|100)+(000|110)\\
	\bW^*&=(111100)=(111|100)+(110|110)+(100|111).
\end{align*}
The tensor products decompose (we indicate the most relevant components by
boxes)
\begin{align*}	
	\bF\otimes\bW&=\boxed{(431|210)}+\boxed{(421|310)} +{(321|320)}\\
	\bE\otimes\bW&={(320|321)}+\boxed{(310|421)}+\boxed{(210|431)} \\
\bF_1\otimes\bW^*&={(542|310)}+{(541|320)}+(532|320)+(531|321)+\boxed{(432|321)}.
\end{align*}
There is no homomorphism from $(310|421)$ to $(431|210)$ (look at the last
but one line in the table in Theorem \ref{class_1sing}), 
hence we can use Proposition \ref{transprop2} (all other assumptions are
clearly satisfied).

Now, all other nonstandard operators in $|1|$-singular case with the lowest
possible weight can be obtained by usual
translations
based on Proposition \ref{transprop1} as follows.
 
\noindent {\bf (A)} $\bE=(210|430)\to\bF=(430|210)\iff 
 \bE'=(210|431)\to \bF'=(431|210)$
\begin{align*}	
 	\bW&=(100000)=(100|000)+(000|100)\\
 	\bW^*&=(111110)=(111|110)+(110|111)\\
 	\bE\otimes\bW&=(310|430)+(210|530) +\boxed{(210|431)}\\
 	\bF\otimes\bW&=(530|210)+\boxed{(431|210)}+(430|310) \\
 	\bE'\otimes\bW^*&=\boxed{(321|541)}+(321|532) +(320|542)\\
 	\bF'\otimes\bW^*&=(543|320)+\boxed{(541|321)}+(532|321).
 \end{align*}
All infinitesimal characters in the last four lines are pairwise different.

\smallskip
 
\noindent {\bf (B)} $\bE=(210|431)\to \bF=(431|210)\iff \bE'=(310|432)\to \bF'=(432|310)$ 
\begin{align*}	
 	\bW&=(100000)=(100|000)+(000|100)\\
 	\bW^*&=(111110)=(111|110)+(110|111)\\
 	\bE\otimes\bW&=(310|431)+(210|531) +\boxed{(210|432)}\\
 	\bF\otimes\bW&=\boxed{(531|210)}+{(432|210)}+(431|310) \\
 	\bE'\otimes\bW^*&=\boxed{(321|542)} +(320|543)\\
 	\bF'\otimes\bW^*&=(543|320)+\boxed{(542|321)}.
 \end{align*}
All infinitesimal characters in the last four lines are pairwise different.
 
\smallskip
 
\noindent {\bf (C)} $\bE=(210|432)\to \bF=(432|210)\iff \bE'=(310|432)\to \bF'=(432|310)$ 
 \begin{align*}	
 	\bW&=(100000)=(100|000)+(000|100)\\
 	\bW^*&=(111110)=(111|110)+(110|111)\\
 	\bE\otimes\bW&=\boxed{(310|432)} +{(210|532)}\\
 	\bF\otimes\bW&={(532|210)}+\boxed{(432|310)} \\
 	\bE'\otimes\bW^*&=(421|542)+{(420|543)} +\boxed{(321|543)}\\
 	\bF'\otimes\bW^*&=(543|320)+\boxed{(543|321)}+{(542|421)}.
 \end{align*}
All infinitesimal characters in the last four lines are pairwise different.
 
\smallskip
 
\noindent {\bf (D)} $\bE=(310|432)\to \bF=(432|310)\iff \bE'=(410|432)\to \bF'=(432|410)$
\begin{align*}	
 	\bW&=(100000)=(100|000)+(000|100)\\
 	\bW^*&=(111110)=(111|110)+(110|111)\\
 	\bE\otimes\bW&=\boxed{(410|432)}+(320|432) +{(310|532)}\\
 	\bF\otimes\bW&={(532|310)}+\boxed{(432|410)}+(432|310) \\
 	\bE'\otimes\bW^*&=\boxed{(521|542)} +(520|543)+(421|543)\\
 	\bF'\otimes\bW^*&=(543|520)+\boxed{(543|421)}+{(542|521)}.
 \end{align*}
All infinitesimal characters in the last four lines are pairwise different.
 
\smallskip

Finally, to reach higher levels of $|1|$-singular infinitesimal character, we can
exploit Proposition \ref{simple_trans}, i.e., twisting with the 
fundamental representations $\mathbb W=(100000)$ or $\mathbb W=(110000)$.

\section{The proof of Theorem \ref{th7}}
 
We start again in the lowest infinitesimal character $(543210)$.  All
standard homomorphisms (included in the BGG complex) are of the first order,
hence we can lift them.
  
To treat  nonstandard ones, we are going to use again 
one-way translations from $|1|$-singular cases. 
We shall construct four homomorphisms of order $4$ and the homomorphism
of order $7$, cf. Theorem \ref{class_regular}.
In view of Theorem \ref{no_lift}, we cannot construct the remaining 
two homomorphisms of order $4$ (involving the first and the last
Verma module in the scheme), while resolving the homomorphism of order $9$ will need much more effort in the separate section below.

\smallskip 
\noindent
{\bf (A)} Translation $V(410|321)\to V(421|310)\implies V(510|432)\to V(532|410)$	
\\ 
Consider the modules
 \begin{align*}
 	\bE&=(410|321);\; \bE_1=(510|432);\; \bE_2=(520|431)\\
 	\bF&=(421|310);\;\bF_1=(532|410);\;\bF_2=(531|420)\\
 	\bW&=(111100)=(111|100)+(110|110)+(100|111)\\	
\bW^*&=(110000)=(110|000)+(100|100)+(000|110).
 \end{align*}
Their tensor products decompose as follows (again we box the most relevant
summands):
\begin{align*}	
 	\bE\otimes\bW&=(521|421)+\boxed{(520|431)}+(421|431)+
\boxed{(510|432)}+(420|432)\\
 	\bF\otimes\bW&=\boxed{(532|410)}+{(532|320)}+
\boxed{(531|420)}+(531|321)+
 	(432|420)+\\
 	&\qquad(432|321)+(521|421)+(431|421)\\ 	 
 	 \bF_1\otimes\bW^*&={(642|410)}+(543|410)+{(632|510)}+
(632|420)+(542|510)+(542|420)+\\
 	 &\qquad{(532|520)}+\boxed{(532|421)}.
 \end{align*}	

Now, $V\bigl((520|431)+(510|432)\bigr)$ splits off from 
$V(\bE\otimes\bW)$ 
(because these two summands are the only summands in  
regular infinitesimal character), and 
$V(\bE_1)=V(510|432)$  embeds to $V(\bE\otimes\bW)$.
Moreover $V[(532|410)+(531|420)]$ splits off from $V(\bF\otimes\bW)$, and
$V(\bF\otimes\bW)$ projects to $V(\bF_1)=V(532|410)$. The condition on the
central action on the summands from Proposition \ref{transprop2} is 
satisfied, too.  
Finally, there is no nontrivial homomorphism from $V(520|431)$ to
$V(532|410)$, see the classification in Theorem \ref{class_regular}. 

Hence, by Proposition \ref{transprop2} 
we have arrived at the nontrivial homomorphism from $V(510|432)$ to $V(532|410)$.
 
\smallskip
 
\noindent
{\bf (B)} Translation $V(310|421)\to V(321|410)\implies V(410|532)\to V(432|510)$

\noindent Consider modules
 \begin{align*}
 	\bE&=(310|421); \; \bE_1=(410|532);\; \bE_2=(420|531)\\
 	\bF&=(321|410);\; \bF_1=(432|510);\;\bF_2=(431|520)\\
 	\bW&=(111100)=(111|100)+(110|110)+(100|111)\\	
\bW^*&=(110000)=(110|000)+(100|100)+(000|110).
 \end{align*}
We have
 \begin{align*}	
 	\bE\otimes\bW&=(421|521)+{(421|431)}+\boxed{(420|531)}+(420|432)+\\
 &\qquad {(321|531)}+{(321|432)}+\boxed{(410|532)}+(320|532)\\
\bF\otimes\bW&=\boxed{(432|510)}+(432|420)+\boxed{(431|520)}+(431|421)+
(421|521)\\
\bF_1\otimes\bW^*&={(542|510)}+{(532|610)}+(532|520)+(432|620)+\boxed{(321|410)}.    
   \end{align*}	
 
There is no nontrivial homomorphism from $V((420|531))$ to $V((432|510))$,
and all the other assumptions in Proposition \ref{transprop2} are valid,
hence we get a nontrivial homomorphism from $V((410|532))$ to
$V((432|510))$.

\smallskip
\noindent
{\bf (C)} Translation $V(320|431)\to V(430|321)\implies V(320|541)\to V(540|321)$
\\
Consider modules
\begin{align*}
 	\bE&=(320|431);\;\bF=(430|321) ;\; \bE_1=(320|541);\; 
\bE_2=(420|531)\\
 	\bF&=(430|321);\;\bF_1=(432|510);\;\bF_2=(431|520)\\
 	\bW&=(110000)=(110|000)+(100|100)+(000|110).
\end{align*}
We have
\begin{align*}	
 \bE\otimes\bW&=(430|431)+{(421|431)}+\boxed{(420|531)}+(420|432)+\\
 	&\qquad {(321|531)}+{(321|432)}+\boxed{(320|541)}+(320|532)\\
 	\bF\otimes\bW&=\boxed{(540|321)}+(531|321)+
\boxed{(530|421)}+(431|421)+(430|431).
 \end{align*}	

There is no nontrivial homomorphism from $V((420|531))$ to $V((540|321))$,
and all other assumption of Proposition \ref{transprop2} are valid,
hence we get a nontrivial homomorphism from $V((320|541))$ to $V(540|321))$.
 
\smallskip
\noindent 
{\bf (D)} Translation $V(321|430)\to V(431|320)\implies V(321|540)\to V(541|320)$
\\
Consider modules
 \begin{align*}
 	\bE&=(321|430); \; \bE_1=(321|540),\; \bE_2=(421|530)\\
 \bF&=(431|320);\;\bF_1=(541|320);\;\bF_2=(531|420)\\
 	\bW&=(110000)=(110|000)+(100|100)+(000|110).
 \end{align*}
We have
 \begin{align*}	
 	\bE\otimes\bW&=(431|430)+\boxed{(421|530)}+(421|431)
 	+ \boxed{(321|540)}+(321|531)\\
 	\bF\otimes\bW&=\boxed{(541|320)}+(532|320)+\boxed{(531|420)}+
(531|321)+
(432|420)+(432|312)+\\&\qquad  (431|430)+(431|421).
 \end{align*}	
 
There is no nontrivial homomorphism from $V((421|530))$ to $V((541|320))$,
and all other assumptions of Proposition \ref{transprop2} are valid,
hence we get a nontrivial homomorphism from $V((321|540))$ to
$V((541|320))$.
 
\smallskip\noindent
{\bf (E)} Translation $V(310|432)\to V(432|310)\implies V(310|542)\to V(542|310)$ (7th order)
\\
Consider modules
 \begin{align*}
 	\bE&=(310|432); \; \bE_1=(310|542);\; \bE_2=(410|532);\;
\bF_1=(542|310)\\
 	\bF&=(432|310);\;\bF_1=(542|310); \;\bF_2=(532|410)\\
 	\bW&=(110000)=(110|000)+(100|100)+(000|110).
 \end{align*}
We have
 \begin{align*}	
 	\bE\otimes\bW&=(420|432)+(321|432)+\boxed{(410|532)}+
(320|532)+\boxed{(310|542)}\\
 	\bF\otimes\bW&=\boxed{(542|310)}+\boxed{(532|410)}+(532|320)+
(432|420)+(432|321).\\
\end{align*}

There is no nontrivial homomorphism from $V((410|532))$ to $V((542|310))$,
and all other assumptions of Proposition \ref{transprop2} are valid,
hence we get a nontrivial homomorphism from $V((310|542)$ to
$V((542|310))$.

\smallskip
 
At the moment, we have verified statement of Theorem \ref{th7} 
for the lowest infinitesimal character.
The proof of the statement for higher level infinitesimal characters is 
finished by standard
translation procedure using the defining representation $\mathbb W=(100000)$ 
and Proposition \ref{transprop1}.

\section{The nonexistence of the longest arrows}

In Theorem \ref{no_lift}, we settled two of the three claims about non-existence of coverings of the singular vectors in the standard Verma modules by reducing the problem to the longest arrow Verma module morphism in the regular pattern in the $(2,2)$-Grassmannian case.

Now, we are going to deal with the longest arrow corresponding the 9th order operator in the $(3,3)$-Grassmannian case directly, mimicking the $(2,2)$-Grassmannian case, which we heavily simplify first. 

\subsection{The $(2,2)$ case---an explicit obstruction certificate}
\label{subsec:gr22-certificate}

	The simplification consists of restricting our attention to the $\mathfrak g_0$ highest weight lines when applying the $\mathfrak g_1$ action. On this very simple case, we shall clarify the general formalism leading to a proof that no cover of a given singular vector of Verma modules can exist. 
	
		\subsection{$\mathfrak g_0$ highest lines}
	The $(2,2)$ Grassmannian geometry comes with the infinitesimal symmetries in the graded Lie algebra $\mathfrak g= \mathfrak g_{-1}\oplus \mathfrak g_0\oplus \mathfrak g_1= \mathfrak{sl}(4,\mathbb R)$.
	The basis of $\mathfrak g$ is given by elementary matrices $E_{ij},1\leq i,j\leq 4.$
	A basis of  $\mathfrak g_{-1}$ will be denoted by
	$$
	Y_{ai}=E_{a+2,i},\qquad a,i=1,2.
	$$

	Let $U$ and $V$ be two-dimensional and identify
	\[
	\mathfrak g_{-1}\cong U\otimes V^*.
	\]
	  Note that if there is a nontrivial semiholonomic singular vector covering the holonomic one, we
can choose it to be  $\mathfrak g_0$-invariant.

	In $\otimes^4 \mathfrak g_{-1}=(\otimes^4 U)\otimes (\otimes^4 V^*),$ the semisimple part of the Levi factor  $\mathfrak g_0$ is realized as 
	  the product of $\mathfrak{sl}(U)$ and $\mathfrak{sl}(V^*),$  so we can discuss the action separately on the left and right factor.

	  Now, Schur-Weyl duality splits $\otimes^4 U$ into three summands but due to invariance of the singular vector, only
the Specht module of shape $(2,2)$ is relevant (see \cite{Ful} for the theory of the representations of the symmetry groups $S_n$ and its link to the representations of $\mathfrak{sl}(U)$).   The dimension $d$ of this space is computed by using the so called hook-rule, i.e., we fill the $n$ boxes in the Young tableau with the numbers $k_i$ of boxes seen on the right and below, and divide $n!$ by the product of all of them. In our $(2,2)$-case, this means 
	$$
	d = \frac{4!}{3\cdot 2\cdot 2 \cdot 1} = 2
	.$$
	As well known, the so called standard Young tableaux are obtained by filling them by letters 1 through $n$, each exactly ones, and keeping the rule that both along the rows and the columns the sequences must grow. In our case, we have two such tableaux, and we can represent them as linear combinations of words in letters 1 through $k$ ($k=2$ in our case), which are the indices indicating which rows the individual letters are in. In our case, there is the space generated by 6 such words. But to encode the proper space of highest-weight lines, we start with the two integral words $1122$ and $1212$ corresponding to the standard Young tableaux and antisymmetrize the columns in them. Thus, we obtain the two word codes that represent the highest weight vectors properly:
	\begin{align*}
		h_1&=1122-1221-2112+2211,\\
		h_2&=1212-1221-2112+2121.
	\end{align*}
	This approach is well-matched to the Schur-Weyl duality, combining the representation of the symmetry group $S_4$ with the tensorial representations of $SL(U)$.
	
	Acting by a $\mathfrak g_0$ highest weight vector $z$ of $\mathfrak g_1$ on a singular vector in our semiholonomic Verma module, we should obtain zero. But hitting a $\frak g_0$-highest weight vector, we always obtain a $\mathfrak g_0$-highest weight vector again. Thus, in our case, the target of the action map will be highest weight vectors of the $(2,1)$ Young tableau type (notice that there is no other option of removing one box from the tableau of type $(2,2)$). The target vector space is thus again 2-dimensional ($3!/(3\cdot 1\cdot 1)$) and is generated by words
	\begin{align*}
		k_1&=112-211,\\
		k_2&=121-211
		.\end{align*}
	
	The general formula for the $\mathfrak g_1$-action is well known, see e.g. \cite{ES}.
	In general, for a $P$-module $\mathbb E$, the elements of $\bigotimes^k{\mathfrak g_{-1}}\otimes{\mathbb E}^*$
	may be written as linear combinations of terms
	$y_1\otimes y_2\otimes\cdots\otimes y_k\otimes e$ for
	$y_1,y_2,\ldots,y_k\in{\frak g_{-1}}$ and $e\in{\Bbb E}^*$.
	By using the Poincar\'e-Birkhoff-Witt procedure, the action of
	$z\in{\frak g}_1$ on such an element is given by
	\begin{equation}\label{z}
		\begin{array}l
			\displaystyle\sum_{1\leq p<q\leq k}y_1\otimes\cdots\otimes\hat{y_p}
			\otimes\cdots\otimes
			[[z,y_p],y_q]\otimes
			\cdots\otimes y_k\otimes e\\[17pt]
			\quad+\displaystyle\sum_{1\leq p\leq k}y_1\otimes\cdots\otimes\hat{y_p}\otimes
			\cdots\otimes y_k\otimes[z,y_p]e\end{array}
		\in\textstyle\bigotimes^{k-1}{\frak g}_-\otimes{\Bbb E}^*,
	\end{equation}
	where the hats over the symbols indicate their omission, while the double commutator appears instead of the omitted $y_q$.
	 
	 A choice of a highest weight vector in $\mathfrak g_1$ depends on a choice of a Borel subalgebra in 
	 $\mathfrak g_0.$ Choosing $E_{21}$ and $E_{34}$  as rising elements. the highest weight vector is $z=E_{24}.$
	 
	The double commutator in the action formula \eqref{z} is then given by,
	\[
	[[z,Y_{ai}],Y_{ck}]
	=
	-\delta_{a2}\delta_{k2}Y_{ci}
	-\delta_{i2}\delta_{c2}Y_{ak}
	.\]

	Since we are dealing with the trivial representation $\mathbb E$, only the first sum survives in the latter formula. 
	We work with the order $k=4$, and elements $y_i\in U\otimes V^*$. All the above considerations are also valid for the 2-dimensional $V^*$, and as a representation space it is naturally isomorphic to $U$. Thus, the space of our $\mathfrak g_0$-highest weight lines is 4-dimensional, generated by tensor products of already discussed highest weight vectors for both $U$ and $V^*$, and we consider them encoded by the words built of the 2-dimensional bases $h_1,\ h_2$.

	For $p\in\{1,2,3,4\}$, let $A_p$ delete the $p$-th letter when that letter is
	$2$, and give zero otherwise.  For $p<q$, let $B_{pq}$ act only when the $q$-th
	letter is $2$: first replace that letter by the original $p$-th letter and
	then delete position $p$ (valued zero otherwise).
	
	As already mentioned, after application of the chosen highest-weight element $z\in\mathfrak
	g_1$, the relevant target multiplicity space for $U$ has shape $(2,1)$, it is again
	two-dimensional, and similarly for $V^*$.
	A straightforward inspection of the action of $z$ on the tensor product of the two multiplicity spaces, as given in \eqref{z}, leads to the following formula:
	
	\begin{lemma} The action of $z$ is 
		represented by the linear mapping $L$,
		\[
		L=-\sum_{p<q}\bigl(A_p\otimes B_{pq}+B_{pq}\otimes A_p\bigr)
		\]
		represented by a 4x4 matrix in our chosen basis, explicitly computed below.
	\end{lemma}
	
	\begin{proof}
		Elements in the space $W= \otimes^4\mathfrak g_{-1}$ can be described  by a product of words using identification
		\begin{equation*}
			Y_{a_1i_1}\otimes	Y_{a_2i_2}\otimes	Y_{a_3i_3}\otimes	Y_{a_4i_4}
			\simeq
			(a_1a_2a_3a_4)\otimes(i_1i_2i_3i_4)
		\end{equation*}
		
		Let us fix a pair $(p,q)$ with $1\leq p<q \leq 4.$
		and consider the action of the first  half 
		$\delta_{a_p,2}\delta_{i_q,2}Y_{a_qi_p}$ of the corresponding summand in the action of the double commutator.  The action
		is nontrivial  iff  $a_p=2$ and $i_q=2$, and if this is true, the action 
		deletes  $Y_{a_pi_p}$  and inserts  $Y_{a_qi_p}$ instead of $Y_{a_qi_q}.$

		The action of the operator  $A_p\otimes B_{pq}$ is also nontrivial  iff  $a_p=2$ and $i_q=2,$
		and if this holds, the action   deletes the letter $a_p$   
		and inserts  $i_p$ instead of $ i_q.$ Hence  both  actions coincide.
		
		The other summand can be treated in the same way.
	\end{proof}

    Restricted to highest weight spaces these maps are
\[
\begin{array}{ll}
 A_1=\begin{pmatrix}-1&-1\\0&1\end{pmatrix},&
 A_2=\begin{pmatrix}0&1\\-1&-1\end{pmatrix},\\[12pt]
 A_3=\begin{pmatrix}1&0\\-1&-1\end{pmatrix},&
 A_4=\begin{pmatrix}1&0\\0&1\end{pmatrix},
\end{array}
\]
and
\[
\begin{array}{lll}
 B_{12}=\begin{pmatrix}0&1\\-1&-1\end{pmatrix},&
 B_{13}=\begin{pmatrix}1&0\\0&1\end{pmatrix},&
 B_{14}=\begin{pmatrix}-1&-1\\1&0\end{pmatrix},\\[12pt]
 B_{23}=\begin{pmatrix}1&0\\-1&-1\end{pmatrix},&
 B_{24}=\begin{pmatrix}0&1\\1&0\end{pmatrix},&
 B_{34}=\begin{pmatrix}1&0\\0&1\end{pmatrix}.
\end{array}
\]

In order to obtain the 4x4 matrix $L$ we have to consider the tensor products. For example, the $p=1$, $q=2$ contribution will combine $A_1$ and $B_{12}$ as 
$$
A_1\otimes B_{12} = \begin{pmatrix} -B_{12} & - B_{12}\\ 0 & B_{12} \end{pmatrix} = \begin{pmatrix} 0&-1&0&-1\\ 1&1&1&1 \\ 0&0&0&1 \\ 0&0&-1&-1\end{pmatrix}
$$
Summing all the 12 contributions (6 options for $p$, $q$ and two terms for each), we arrive at
\[
 L=
 \begin{pmatrix}
 -2&-1&-1&-2\\
  2& 1& 1& 2\\
  2& 1& 1& 2\\
  0& 0& 0& 0
 \end{pmatrix}.
\]
Thus, the action of $z$ on the highest weight $x = x_{11}h_1\otimes h_1 + x_{12}h_1\otimes h_2 + x_{21}h_2\otimes h_1 + x_{22}h_2\otimes h_2$ is given by multiplying $L\cdot x$, where $x$ is the vector of the latter coefficients. 

Next we want to understand which linear combination of words is corresponding to the holonomic symmetric weight vector, represented by the square of the determinant. We can see it with the help of an invariant scalar product on the space of words. Indeed, we declare the distinct simple words orthonormal. This scalar product is clearly invariant with respect to the symmetric group $S_4$ and the relevant Gram matrix of this symmetric bilinear form in our basis $h_1$ and $h_2$ is 
\[
 G=
 \begin{pmatrix}
 4&2\\2&4
 \end{pmatrix}.
\]
Passing to the holonomic case, there is only one-dimensional space of invariant linear functionals $s$ left and thus we may write this as $s(x\otimes y) = \langle x,y\rangle_G$, which is clearly the invariant choice. The holonomic symbol functional is thus
\[
 s=\operatorname{vec}(G)
   =(4\ 2\ 2\ 4),
\]
where $\operatorname{vec}$ denotes row-major vectorization.  Notice that the
canonical diagonal $S_4$-invariant tensor is represented by $G^{-1}$, and
$s(\operatorname{vec}(G^{-1}))=2$, the dimension of the Specht module
$[2,2]$.

Now, in order to get a certificate of non-existence of a $\mathfrak p$-highest weight vector, we need any linear form $\ell$ such that 
\begin{equation}\label{certificate}
\ell\cdot L = s
.\end{equation}
Indeed, then for any vector $v$ in the kernel of $L$ its holonomic projection $s\cdot v$ vanishes.
Now, there is the obvious solution to \eqref{certificate}:
\[
 \ell=(-2\ 0\ 0\ 0)
.\]
Consequently, if $v$ is a semiholonomic singular candidate, so that $Lv=0$,
then
\[
 s(v)=\ell\cdot L\cdot v=0.
\]
Hence no vector in $\ker L$ can project nontrivially to the one-dimensional
holonomic symbol line generated by $(\det Y)^2$.

If we want to stay in the integer matrix computation, it is convenient to rewrite
the same certificate in $2\times2$ block form.  We write $E$ for the matrix form of our $\ell$,
\[
 E=\begin{pmatrix}-2&0\\0&0\end{pmatrix}.
\]
and then $\ell L=s$ is equivalent to
\begin{equation}
 -\sum_{p<q}
 \left(A_p^{\mathsf T}E B_{pq}
       +B_{pq}^{\mathsf T}E A_p\right)
 =G.
 \label{eq:gr22-block-certificate}
\end{equation}
Indeed, the six summands on the left-hand side are respectively
\begin{align*}
 C_{12}&=\begin{pmatrix}0&-2\\-2&-4\end{pmatrix},&
 C_{13}&=\begin{pmatrix}-4&-2\\-2&0\end{pmatrix},&
 C_{14}&=\begin{pmatrix}4&4\\4&4\end{pmatrix},\\[8pt]
 C_{23}&=\begin{pmatrix}0&2\\2&0\end{pmatrix},&
 C_{24}&=\begin{pmatrix}0&0\\0&4\end{pmatrix},&
 C_{34}&=\begin{pmatrix}4&0\\0&0\end{pmatrix},
\end{align*}
where
\[
 C_{pq}=-\left(A_p^{\mathsf T}E B_{pq}
                   +B_{pq}^{\mathsf T}E A_p\right).
\]
Adding them gives
\[
 \sum_{p<q}C_{pq}
 =\begin{pmatrix}4&2\\2&4\end{pmatrix}
 =G,
\]
which proves \eqref{eq:gr22-block-certificate} by direct integer arithmetic.
Thus the non-lifting obstruction in the $(2,2)$ case is already visible in a
single nonzero coordinate of the dual certificate.  

The $(3,3)$ calculation
uses exactly the same identity, with the two-dimensional multiplicity spaces
replaced by the much higher dimensional Specht module spaces. Let us look at this case in more detail, though the computation itself is left to computer program code.

\subsection{The $(3,3)$ case---an explicit obstruction certificate}
\label{subsec:gr33-certificate}

The same procedure applies to the $(3,3)$--Grassmannian case, but the dimensions increase considerably. We therefore only summarize the computation and refer to the ancillary code files for its exact implementation and verification.

Let $U$ and $V$ be three-dimensional and write
\[
\mathfrak g_{-1}=U\otimes V^*.
\]
The flat operator under consideration has order $9$ and its holonomic principal symbol is, up to scale,
\[
(\det Y)^3\in S^9(U\otimes V^*).
\]
For either of the two tensor factors, the relevant weight space in the ninth tensor power consists of words containing three copies of each of the letters $1,2,3$:
\[
C_9(3,3,3).
\]
Its dimension is already
\[
\dim C_9(3,3,3)
 =\frac{9!}{3!3!3!}=1680.
\]
The subspace of the $\mathfrak g_0$ highest vector lines is the highest-weight multiplicity space
\[
H_{333}\simeq [3,3,3].
\]
By the hook-length formula,
\[
\dim H_{333}=42.
\]
Equivalently, Schur--Weyl duality identifies this space with the multiplicity space of the one-dimensional $GL(3)$-module
\[
\mathbb S_{(3,3,3)}U=(\det U)^3.
\]
Thus, for the two factors $U$ and $V^*$, the space of possible $\mathfrak g_0$-singular semiholonomic covers has dimension
\[
\dim(H_{333}\otimes H_{333})=42^2=1764.
\]

We choose the highest-weight vector $z=E_{36}\in\mathfrak g_1$. Its action lowers the tensor degree by one. In one tensor factor the relevant target weight is $(3,3,2)$; the corresponding full word space has dimension
\[
\dim C_8(3,3,2)
  =\frac{8!}{3!3!2!}=560.
\]
The highest-weight subspace is
\[
H_{332}\simeq[3,3,2],
\qquad
\dim H_{332}=42.
\]
The equality of the two dimensions may also be seen directly from the symmetric-group branching rule: the rectangular Young diagram $(3,3,3)$ has a unique removable corner.

Exactly as in the $(2,2)$ example, the double commutator
\[
[[z,Y_{ai}],Y_{ck}]
 =
-\delta_{a3}\delta_{k3}Y_{ci}
-\delta_{i3}\delta_{c3}Y_{ak}
\]
gives two elementary operations on the word spaces. For each tensor index position $p$ we obtain a map
\[
A_p:H_{333}\longrightarrow H_{332},
\]
and for each $p<q$ a map
\[
B_{pq}:H_{333}\longrightarrow H_{332}.
\]
After choosing the right word bases (by taking the standard tableaux and antisymmetrizing the columns), all $A_p$ and $B_{pq}$ are $42\times42$ integer matrices. There are nine matrices $A_p$ and $\binom92=36$ matrices $B_{pq}$. The action of $z$ on the full multiplicity space is therefore represented by the $1764\times1764$ matrix
\[
L=
-\sum_{p<q}
\bigl(
 A_p\otimes B_{pq}+B_{pq}\otimes A_p
\bigr).
\]
This is exactly the same formula as in the $(2,2)$ computation; only the dimensions have changed.

Let $G$ denote the $42\times42$ Gram matrix of our integral polytabloid basis of $[3,3,3]$. The holonomic symmetrization is represented, as before, by the invariant functional
\[
s=\operatorname{vec}(G).
\]
Thus a semiholonomic lift of the nonzero holonomic symbol could exist only if there were a vector $v$ satisfying
\[
Lv=0,\qquad s(v)\neq0.
\]
To exclude this, it suffices to exhibit a dual certificate showing that $s$ belongs to the row space of $L$.

The computation produces such a certificate exactly. After clearing denominators, its $1764$ coefficients may be reshaped into an integer $42\times42$ matrix $E$, solving the system of linear equations. Of its $1764$ entries, $567$ are nonzero. The resulting identity is
\[
-\sum_{p<q}
\left(
 A_p^{\mathsf T}E B_{pq}
 +
 B_{pq}^{\mathsf T}E A_p
\right)
=
2G.
\]
Equivalently, in vectorized form (in the program code, all vectors are represented as columns),
\[
L^{\mathsf T}\ell=2s,
\qquad
\operatorname{vec}(E)=\ell.
\]
Consequently,
\[
Lv=0
\quad\Longrightarrow\quad
s(v)=0,
\]
and hence no semiholonomic singular vector can project to the nonzero holonomic symbol $(\det Y)^3$.

We have created python code which checks all ingredients in this argument by exact integer arithmetic. 

The code \texttt{grassmannian\_lift\_step1.py} constructs the polytabloid bases, verifies the highest-weight conditions, computes the Gram matrix and the $A_p,B_{pq}$ matrices, and checks their action directly in the full word spaces. 

The code \texttt{grassmannian\_lift\_step2.py} assembles the operator $L$ and verifies the displayed certificate both in its $42\times42$ block form and by the corresponding full Kronecker-product computation. A modular computation gives rank $1089$ for $L$ over several good primes; this is useful diagnostic information but is not needed for the proof. The non-lifting statement follows solely from the exact certificate above.

In the standard Google Colab environment, the above mentioned ancillary files have to be uploaded to the temporary file repository, and then the following sequence of executed commands provides the quoted results:

\begin{verbatim}
!pip install sympy numpy scipy numba
!python grassmannian_lift_step1.py \
    --d 3 --m 3 \
    --output gr33_step1_matrices.json.gz

Source shape (3, 3, 3): 42 standard tableaux; word dimension 1680
Target shape (3, 3, 2): 42 standard tableaux; word dimension 560
PASS: all basis vectors annihilated by all simple raising operators
PASS: exact inverses of both 42- and 42-dimensional Gram matrices
PASS: symbol G and diagonal S_9-invariant tensor G^-1; symbol normalization = 42
PASS: independent matrix-unit double-commutator agrees with A/B maps
PASS: 9 deletion matrices A_p reconstructed exactly
PASS: 36 replacement matrices B_pq reconstructed exactly
PASS: all 42-by-42 matrix entries integral
Saved exact matrices and Gram matrices: gr33_step1_matrices.json.gz
PASS: all first-stage checks; elapsed 12.82 s

!python grassmannian_lift_step2.py \
  --input gr33_step1_matrices.json.gz \
  --certificate reproduced_certificate.json

Integer L: (1764, 1764), nonzero entries 506627; symbol nonzeros 1668
Solving L.T ell = s over GF(32003)...
Modular rank(L)=1089; augmented system inconsistent=False
EXACT non-lifting certificate verified: {'type': 'obstruction', 'scale': 2,
   'nonzero_coefficients': 567, 'identity': 'L.T (entries / scale) = vec(G)',
   'all_residuals_zero': True}
Certificate saved: reproduced_certificate.json
PASS: independent full-matrix exact check (506627 nonzeros)
PASS: conclusive exact obstruction certificate; elapsed 17.29s

!python grassmannian_lift_step2.py \
    --verify-only \
    --input gr33_step1_matrices.json.gz \
    --certificate reproduced_certificate.json

PASS: independently verified exact certificate: {'type': 'obstruction',
   'scale': 2, 'nonzero_coefficients': 567,
   'identity': 'L.T (entries / scale) = vec(G)', 
   'all_residuals_zero': True}
PASS: independent exact full-matrix cross-check (506627 nonzeros)
\end{verbatim}

The important point is that the found certificate $\ell$ is fully checked in full integer arithmetic, although it was computed in the modular one. Thus, there is no doubt that the found linear form $\ell$ is a genuine certificate of the non-existence of a covering of the holonomic singular vector.

\subsection{Concluding remarks}
The python code verification of the existence of the certificate has been supported by discussion with the chatGPT, version Plus. Of course, it is very tempting to look for the general dimensions $(n,n)$-Grassmannian geometries. But there is no hope to continue the same way: already for $n=4$, we face the multiplicity 24\,024 of the $\mathfrak g_0$ highest weight vector lines in the semiholonomic Verma modules in the relevant order 16. Thus the size of the matrix $L$ would become 577\,152\,576 x 577\,152\,576. Still, we hope that better understanding of the structural essence of our argumentation could lead to a breakthrough in the future. 

Even more interesting question stays open here: in the conformal geometry de Rham sequences, the longest arrow always exists as a truly exceptional invariant differential operator. This is caused by the fact, that the universal differential operator into the semiholonomic jets is not surjective and, thus, it is enough the check the $P$-equivariance of our formulae on the image (a subspace defined by Bianchi and Ricci identities and all their differential consequences). 

In the conformal Riemannian geometry, the holonomic formula demonstrates such a behaviour, at least for the longest arrow in the de Rahm sequence (and it is known that higher then the critical powers of the Laplacian do not exist). Similar behavior in the almost quaternionic case was found in the work of N\'avrat, \cite{Navrat}, but there the natural torsion had to be involved in modification of the flat geometry operator formula. All the interesting nonstandard operators in the de Rham sequence are of fourth order there. 

Actually, we have been surprised, that in the $(3,3)$ case four of the six similar fourth order operators actually allow the semiholonomic module homomorphism cover!

Although a similar manipulation with the torsion correction terms as in N\'avrat's work could be possible, this seams to be a challenging programme.

\subsection{Declaration of generative AI and AI-assisted technologies in the manuscript preparation process}
During the preparation of this work the authors used chatGPT 5 version Plus to double check the technique applied in this very last section, and prepare the Python code for final computational check, as described in detail above. After using this tool, the authors reviewed and edited the content as needed and take full responsibility for the content of the published article.

\end{document}